\documentclass[11pt]{article}
\usepackage[a4paper,left=2cm,right=2cm,top=3cm,bottom=3cm]{geometry}

\usepackage{xcolor}
\usepackage{url} % 
\usepackage{caption}
\usepackage{amsmath}
\usepackage{amssymb}
\usepackage{mathtools}
\usepackage{amsthm}
\usepackage{microtype}
\usepackage{graphicx}
\usepackage{subfigure}
\usepackage{booktabs} % for professional tables
\usepackage{bbm}
\usepackage{tabularray}
\usepackage{algorithm}% http://ctan.org/pkg/algorithms
\usepackage{algpseudocode}% http://ctan.org/pkg/algorithmicx
\usepackage{authblk}
\title{Zigzag Persistence of Neural Responses to Time-Varying Stimuli}

\author[1,2]{Yuri Gardinazzi}
\author[1]{Alessio Ansuini}
\author[3]{Eugenio Piasini}
\author[2]{Fabio Anselmi}
\author[1]{Matteo Biagetti\thanks{For correspondence: matteo.biagetti@areasciencepark.it}}
\affil[1]{Area Science Park, Padriciano 99, Trieste}
\affil[2]{Department of Mathematics, Informatics and Geosciences\\ 

University of Trieste, Via Economo 12/3, Trieste}
\affil[3]{SISSA, Via Bonomea 265, Trieste }

\definecolor{ForestGreen}{rgb}{0.13, 0.55, 0.13}
\definecolor{airforceblue}{rgb}{0.36, 0.54, 0.66}
\definecolor{orange}{rgb}{1.0, 0.5, 0.0}
\definecolor{amethyst}{rgb}{0.6, 0.4, 0.8}
\definecolor{awesome}{rgb}{1.0, 0.13, 0.32}
\definecolor{chromeyellow}{rgb}{1.0, 0.65, 0.0}

\date{}

\begin{document}

\maketitle

\begin{abstract}
We use topological data analysis to study neural population activity in the Sensorium 2023 dataset, which records responses from thousands of mouse visual cortex neurons to diverse video stimuli. For each video, we build frame-by-frame cubical complexes from neuronal activity and apply zigzag persistent homology to capture how topological structure evolves over time. These dynamics are summarized with persistence landscapes, providing a compact vectorized representation of temporal features. We focus on one-dimensional topological features—loops in the data—that reflect coordinated, cyclical patterns of neural co-activation. To test their informativeness, we compare repeated trials of different videos by clustering their resulting topological neural representations. Our results show that these topological descriptors reliably distinguish neural responses to distinct stimuli. This work highlights a connection between evolving neuronal activity and interpretable topological signatures, advancing the use of topological data analysis for uncovering neural coding in complex dynamical systems.
\end{abstract}

\section{Introduction}
A fundamental challenge in neuroscience is understanding how information about the world and an animal's internal state and goals is encoded in the brain. Modern electrophysiological and imaging recording techniques allow for simultaneous recording of very large populations of neurons with single-cell resolution, and have driven the design of methods for the analysis of collective patterns of neuronal activity and their representational properties \cite{saxenaNeuralPopulationDoctrine2019}.

%A fundamental challenge in neuroscience is understanding how the brain represents information across space and time form sensory signals to high-level cognitive tasks. To this aim, researchers, instead of focusing on single neurons, commonly examine patterns within neural population activity and study how different brain regions interact, using both activity-based and connectomic approaches~\cite{sporns2005human,saxenaNeuralPopulationDoctrine2019,bassett2017network}.

A central concept in this research is that of the \emph{neural manifold} -- a low-dimensional geometric embedding of population activity within a high-dimensional neural state space~\cite{gallego2017neural}.
Experimental and theoretical studies suggest that neuronal activity typically lies on such manifolds, whose structure reflects the underlying task variables, e.g. object identity in the case of object recognition in visual cortex~\cite{dicarloHowDoesBrain2012,chungNeuralPopulationGeometry2021,perichNeuralManifoldView2025}.
Neural manifold learning methods---including both linear (e.g., PCA) and non-linear (e.g., Isomap, UMAP) techniques---have been applied to provide insight into systems as diverse as the head direction circuits in  mouse~\cite{chaudhuri2019intrinsic} and zebrafish~\cite{petruccoNeuralDynamicsArchitecture2023}, grid cells~\cite{gardnerToroidalTopologyPopulation2022}, motor regions across multiple species~\cite{gao2017theory,fortunatoNonlinearManifoldsUnderlie2024} and head ganglia during locomotion in C.\ Elegans \cite{brennanQuantitativeModelConserved2019}.

%across hippocampus, motor cortex, and prefrontal cortex, providing insight into neural dynamics in health and disease~\cite{chaudhuri2019intrinsic,gao2017theory}. Moreover, work on head-direction systems, grid cells, and representations in V1 visual cortex further grounded the neural manifold framework in concrete physiological circuits,~\cite{seelig2015neural,fiete2008grid}.

Beyond a geometric characterization of the neural manifold, \emph{topological data analysis} provides another layer of understanding by focusing on shape invariants rather than coordinates.
Persistent homology---a key technique in topological data analysis---captures multi-scale topological features (such as loops and voids) in a robust, noise-tolerant manner \cite{edelsbrunner2000topological,carlsson2009topology}.
Several works have used these tools in the context of computational neuroscience (e.g. \cite{giusti2015clique,riemanncliques}).
However, classical persistent homology assumes a \textit{static} point cloud filtration.
Neuronal activity, by contrast, evolves dynamically, especially under time-varying stimuli or behavior. \emph{Zigzag persistent homology} \cite{carlsson2009zigzag} extends persistence to sequences that allow additions and deletions of points, in our context the neuronal activity in time, enabling the direct tracking of topological structure as it changes. 

In this work, we leverage zigzag persistence to analyze population responses from the Sensorium 2023 dataset \cite{sensorium2023}, which contains records of visual cortex neuronal activity of mouse exposed to grayscale videos. 
We develop a pipeline that takes the neuronal activity for each video, frame-by-frame, and transforms it into a sequence of evolving topological shapes underlying the neuronal activity. We then use zigzag persistence to track how these shapes change over time, focusing on cyclical patterns of neural co-activation. We summarize these evolving features using persistence landscapes~\cite{bubenik2015statistical}, providing a compact, vectorized representation of the topological changes in time.

Our experiments address three questions: (A) do the zigzag descriptors separate repeated presentations of different videos to the same mouse? (B) do they carry signal about coarse video families (e.g., Naturalistic, Gaussian, Waves)? and (C) do they carry mouse identity? Our results show that:

\begin{itemize}
\item[(A1)]Repeated presentations of distinct videos reliably cluster by individual movie, with high Adjusted Rand Index \cite{Hubert1985ComparingP} across video types. 
%This indicates that time-resolved loop structure ($H_1$ zigzag) captures stimulus-specific spatiotemporal organization in the neural responses.

\item[(A2)] %Controls: two perturbations confirm the dependence on spatiotemporal structure. 
Shuffling frame order within each video collapses clustering performance, showing that temporal order is critical. Scrambling the spatial grid (per plane, per video) degrades performance toward chance, showing that spatial contiguity is also critical.

\item[(B)] Linear classifiers trained on our descriptors to predict video types achieve accuracy modestly above chance. 
%This suggests that our descriptor is more sensitive to movie-specific spatiotemporal structure than to coarse stimulus categories.

\item[(C)] Linear classifiers trained on our descriptors for mouse identification perform only slightly above chance.
%, suggesting limited sensitivity to idiosyncratic per-animal features. 
%Because the dataset lacks shared movies across mice, we cannot directly assess cross-mouse generalization for the same stimulus.
\end{itemize}

Overall, these results indicate that applying zigzag persistence to frame-wise activity fields produces compact, interpretable topological signatures that are informative about specific dynamic stimuli.
%Our contributions are:
%(i) a practical pipeline adapting zigzag persistence to frame-by-frame population imaging via 2D superlevel-set constructions and a simplicial adapter, summarized with persistence landscapes;
%(ii) empirical evidence that $H_1$ zigzag descriptors reliably distinguish repeated presentations of individual videos within mouse, while only modestly separating broad stimulus families;
%(iii) control analyses showing that both temporal order and spatial contiguity are necessary for the observed clustering performance.
This approach opens a promising direction for future research, that of
extending topological data analysis to capture the temporal structure of neural representation, potentially across behavioral contexts, brain regions, or in the presence of learning and adaptation in time.

\section{Dataset}\label{sec:dataset}

The Sensorium 2023 movie dataset~\cite{sensorium2023} contains two-photon calcium imaging data recorded from the visual cortex of behaving mice during repeated presentations of video clips. For each mouse, neuronal activity was recorded simultaneously from multiple imaging planes (z-slices) in V1. The dataset provides, for each video frame, single-trial fluorescence responses for all recorded neurons in each plane, along with 2D spatial coordinates of cell bodies and behavioural data (e.g. pupil diameter, treadmill velocity).

We consider a subset of $5$ mice in the training set release. We select videos with the following classes: \textsc{naturalistic}, \textsc{gaussian}, \textsc{waves} and \textsc{moving dot}. Each class features different videos with several repeats.\footnote{For more details, see Appendix \ref{app:dataset}.}
%\footnote{Repeats are useful to compare neuronal activity of a given mouse with the same baseline visualization. Notably, we were unable to find same-video, different-mouse cases, i.e., each mouse watched different videos from the others.} 
 In this subset, each scan comprises $10$ $z$-planes, and each plane contains on average approximately $800$ neurons, yielding several thousand units per mouse. Stimuli are grayscale videos presented at a fixed frame rate; responses are provided preprocessed by the authors of the dataset, aligned to stimulus frame times.\footnote{We do not include behavioural data in our analysis, which would allow for a more fine-grained examination of neural responses. We leave this for future work.}

We organize the data as follows: for each mouse $a$, plane $p$, video $v$, and frame $t$, we extract the vector $r_{a,p,v}(t)$ of single-trial neural responses across neurons in that plane. We keep only frames belonging to the active video (excluding inter-trial blanks). For downstream analysis, we transform $r$ into $\delta_{a,p,v}(t) = (r_{a,p,v}(t) - \bar{r}_{a,p,v})/\bar{r}_{a,p,v}$,
where the average is across frames, such that $\delta$ is a normalized response, which ranges from $-1$ (no response) to any positive value, and it is zero at average response per mouse. Each $\delta$ is interpolated on a 2D grid using Piecewise Cubic Spline interpolation. These per-frame grids are the inputs to the cubical complexes used by our zigzag pipeline, explained below.

\section{Methodology}\label{sec:method}

We build a sequence of complexes from per-frame 2D activity grids (one per $z$-plane) and extract their evolving topological features via zigzag persistence. 

For each mouse $a$, plane $p$, video $v$, and frame $t$, we have a 2D activity $\delta_{a,p,v}(t)$ on a fixed 2D grid (Sec.\ref{sec:dataset}).
We form a superlevel-set cubical complex $\mathcal{K}_{a,p,v}(t)$ by selecting all grid points (and their incident edges/squares) whose value exceeds the threshold $\delta_{{\rm threshold}}=0$ (activity above the per-neuron mean across frames). This yields one 2D cubical complex per frame.\footnote{Cubical complexes have also been used in time-varying functional Magnetic Resonance Imaging (fMRI), see \cite{rieckcubical}. }
%\paragraph{From cubical to simplicial complexes.}

For practical use of publicly available code computing zigzag persistence\footnote{We use a combination of \textsc{Dionysus2} \cite{dionysus} and \textsc{fastzigzag} \cite{dey2022fastcomputationzigzagpersistence}, see details in Appendix \ref{app:algo}.}, we map each $\mathcal{K}_{a,p,v}(t)$ to a simplicial complex by closing the corner set of every active square into an abstract simplex and then taking the 2-skeleton:
\begin{equation}
\mathcal{S}_{a,p,v}(t)\;=\;\mathrm{Sk}_2\!\left(\bigcup\{C\ \texttt{active (square/edge/vertex) at } t\}\ \langle V(C)\rangle,\right).
\end{equation}
Here, $C$ ranges over active cubical cells (vertices, edges, squares), $V(C)$ is the set of its grid-corner vertices, and $\langle V(C)\rangle$ denotes the abstract simplex on $V(C)$ (in a simplex tree, inserting $\langle V(C)\rangle$ automatically inserts all its faces). Concretely, we:
(i) insert all vertices (0-simplices);
(ii) insert all active grid edges (1-simplices);
(iii) for each active square with corners ${v_{00},v_{10},v_{01},v_{11}}$, insert the 3-simplex $[v_{00},v_{10},v_{01},v_{11}]$ and retain its 2-skeleton.
This is not a triangulation of the square; it is an abstract “closure” that adds both diagonals and all triangular faces among the four corners.\footnote{For details about the algorithm and a discussion about this choice, we refer to Appendix \ref{app:algo}.}
%\paragraph{Time zigzag via interleaved intersections.}

Having set up the simplicial complexes for each frame, we now construct an interleaved sequence over consecutive frames using intersections at intermediate steps:
\begin{equation}
\mathcal S_{a,p,v}(t_1)\ \hookleftarrow \mathcal S_{a,p,v}(t_1)\,\cap\, \mathcal S_{a,p,v}(t_2)\ \hookrightarrow\ \mathcal S_{a,p,v}(t_2)\ \hookleftarrow\ \cdots
\end{equation}
This “intersection–zigzag” matches our implementation, where even layers contain the common simplices of adjacent frames. We compute zigzag persistent homology in dimension 1 ($H_1$) on this sequence using code based on the implementation of \cite{gardinazzipersistent} for neural network layers.
%\paragraph{ Vectorization with persistence landscapes.}
Each $H_1$ zigzag barcode is mapped to persistence landscapes $\Lambda_k$, where $k$ is the landscape layer, over the time index. We sample the landscapes at a resolution of 50 and retain the first 5 landscape layers ($k=1,\dots,5$), yielding a 250-dimensional vector per plane.

\section{Experiments}

We evaluate whether zigzag descriptors capture stimulus information across three settings: (A) clustering of repeated presentations of distinct videos of the same video type for a given mouse; (B)  classification across video type for a given mouse; and (C) across-mouse classification. Each trial is represented by a 2500-dimensional vector obtained by concatenating the $10$ $z$-planes.
\paragraph{Clustering protocol (A).} For a given mouse and a given video type, we cluster repeated presentations of different videos of that type. Before clustering, we normalise all data vectors in the range $[0,1]$ and reduce the dimensionality to $10$ dimensions via PCA fit on the pooled set under evaluation. We then apply agglomerative clustering (Ward linkage, Euclidean distance) with the number of clusters set to the ground truth. Performance is reported as Adjusted Rand Index, for which chance corresponds to $\approx 0$), mean $\pm$ std over $20$ resamplings. In each resampling, we draw balanced subsets (typically $10$ samples per class when available). As a consistency check, we run two different tests: (i) we permute frames within each video identically across planes (disrupt temporal order, preserve per-frame spatial structure) and (ii) we permute pixel values on each plane once per video (disrupt spatial contiguity, largely preserve marginal activity statistics). Results for these experiments are reported in Table~\ref{tab:exp_ari_main}.\footnote{In Appendix \ref{app:results} we quote an extended range of results for multiple mice and including other metrics (accuracy of cluster identification and adjusted mutual information).}
\begin{table}[h]
\centering
\small
\begin{tabular}{lccc}
\toprule
Condition & Baseline ARI & Frame-shuffled ARI & Grid-scrambled ARI \\
\midrule
  Naturalistic & 0.945 $\pm$ 0.107 & 0.194 $\pm$ 0.118 & 0.180 $\pm$ 0.153 \\
  Gaussian & 0.702 $\pm$ 0.156 & 0.018 $\pm$ 0.024 & 0.081 $\pm$ 0.072 \\
  Waves & 0.619 $\pm$ 0.197 & -0.005 $\pm$ 0.011 & 0.028 $\pm$ 0.004 \\
\bottomrule
\end{tabular}
\caption{Clustering ARI (mean $\pm$ std over 20 runs) for the mouse labelled \textsf{2--10}. Baseline uses intact data. Frame-shuffled permutes frames within each video. Grid-scrambled permutes pixels per plane once per video.}
\label{tab:exp_ari_main}
\end{table}
\vspace{-2mm}
\paragraph{Classification protocol (B) video type.} For a given mouse, we train a linear logistic regression classifier to predict video type (e.g., Naturalistic, Gaussian, Waves) from the 2500-dimensional descriptors. We create stratified 80/20 train/test splits and repeat this procedure over $5$ random splits; we find a  cross-validation accuracy of $\mathbf{0.69 \pm 0.06}$. We also report the f1-score per class and confusion matrix aggregated on test-set predictions across splits (see Figure \ref{fig:classifier} in Appendix \ref{app:results}).
\paragraph{Classification protocol (C) mouse identity. }
We similarly train a linear logistic regression classifier to predict mouse identity from the 2500-dimensional descriptors, pooling all available videos of the selected types across $5$ mice. Because the same video is never shared across mice (due to dataset construction), this task probes whether descriptors predominantly reflect individual mouse features. In this case, the cross-validation accuracy is $\mathbf{0.36 \pm 0.03}$. The aggregated confusion matrix is shown in Figure \ref{fig:classifier} in Appendix \ref{app:results}).
% \begin{center}
% \captionsetup{type=table}
% \label{tab:classifier_results}
% \small
% \begin{tabular}{@{}l l@{}}
% \textbf{Video family} 0.69 ± 0.06 & \textbf{Mouse identity} 0.36 ± 0.03
% \end{tabular}
% \captionof{table}{Cross-validated accuracy for classification protocols.}
% \end{center}
% \begin{table}[ht]
% \centering
% \begin{tabular}{lc}
% \hline
% \textbf{Experiment} & \textbf{Value} \\
% \hline
% Video Family & 0.69 $\pm$ 0.06 \\
% \hline
% Mouse Identity & 0.36 $\pm$ 0.03 \\
% \hline
% \end{tabular}
% \caption{Cross-validated accuracy for classification protocols. Video family results are quoted for mouse \textsc{2-10}.}
% \label{tab:classifier_results}
% \end{table}
\section{Results and Conclusions}

We presented a pipeline that applies zigzag persistent homology to frame-by-frame neural population activity. On Sensorium 2023, $H_1$ zigzag descriptors based on persistence landscapes robustly separate repeated presentations of individual videos within mouse (high ARI; $\approx 0.94$ Naturalistic, $\approx 0.70$ Gaussian, $\approx 0.62$ Waves) and only modestly separate stimulus type ($\approx 70\%$ accuracy). Two controls indicate that temporal order and spatial contiguity are essential: shuffling frames collapses clustering, and scrambling spatial grids drives performance toward chance. Mouse identification is just above chance when classifying mice that viewed disjoint movies ($\approx 35\%$ accuracy), suggesting limited sensitivity to mouse-specific characteristics.

Methodologically, we converted cubical complexes into simplicial complexes by closing active squares with an abstract simplex and retaining its 2-skeleton; this may over-connect cells and accelerate cycle filling. Future work will compare this choice with triangulated and cubical alternatives. The analysis can also be improved by incorporating behavioral data to generally reduce unexplained variance in neural responses. Overall, our results indicate that time-resolved cyclic coordination patterns in neural populations carry stimulus-specific information that zigzag persistence can extract and quantify.
\section*{Acknowledgements}

M.B. thanks Mathieu Carrière and Karthik Viswanathan for useful feedback on the zigzag implementation. M.B., Y.G. and F.A. thank Nikos Karantzas for discussion about the dataset in the early stages of the project. Y.G. is supported by the
Programma Nazionale della Ricerca (PNR) grant J95F21002830001 with title “FAIR-by-design”. E.P. was partly supported by the PRIN Project no. 2022XE8X9E (CUP:G53D23004590001) and by the SISSA 5x1000 IRPEF 2023 funds.

We thank Area Science Park supercomputing platform ORFEO made
available for conducting the research reported in this paper and the technical support of the Laboratory of Data Engineering staff.

\bibliographystyle{utphys}
\bibliography{gtml2025_workshop}

\appendix

\section{More info on Sensorium (2023)}\label{app:dataset}

\begin{figure}[h!]
    \centering
    \includegraphics[width=0.45\textwidth]{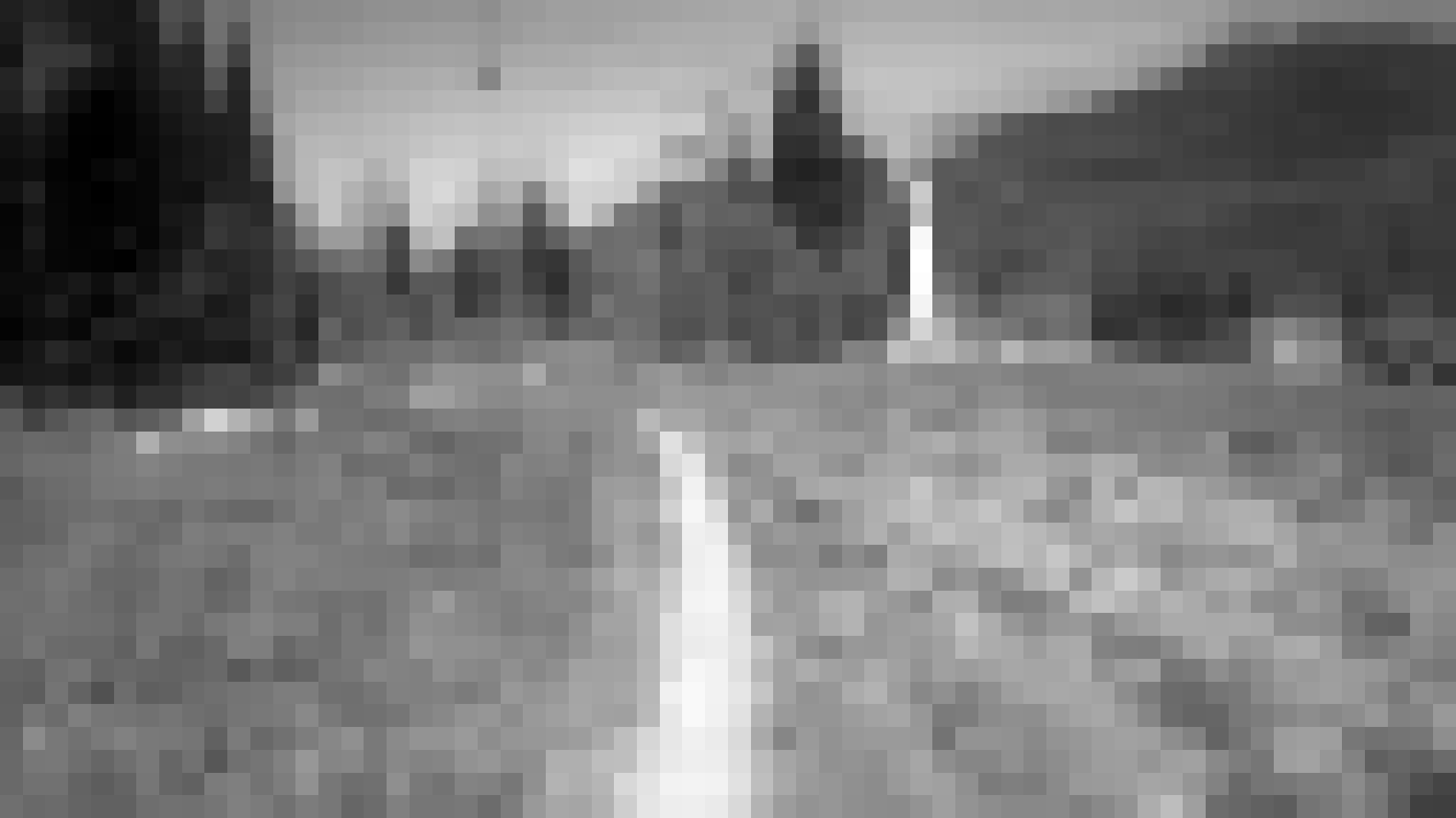}
    \includegraphics[width=0.45\textwidth]{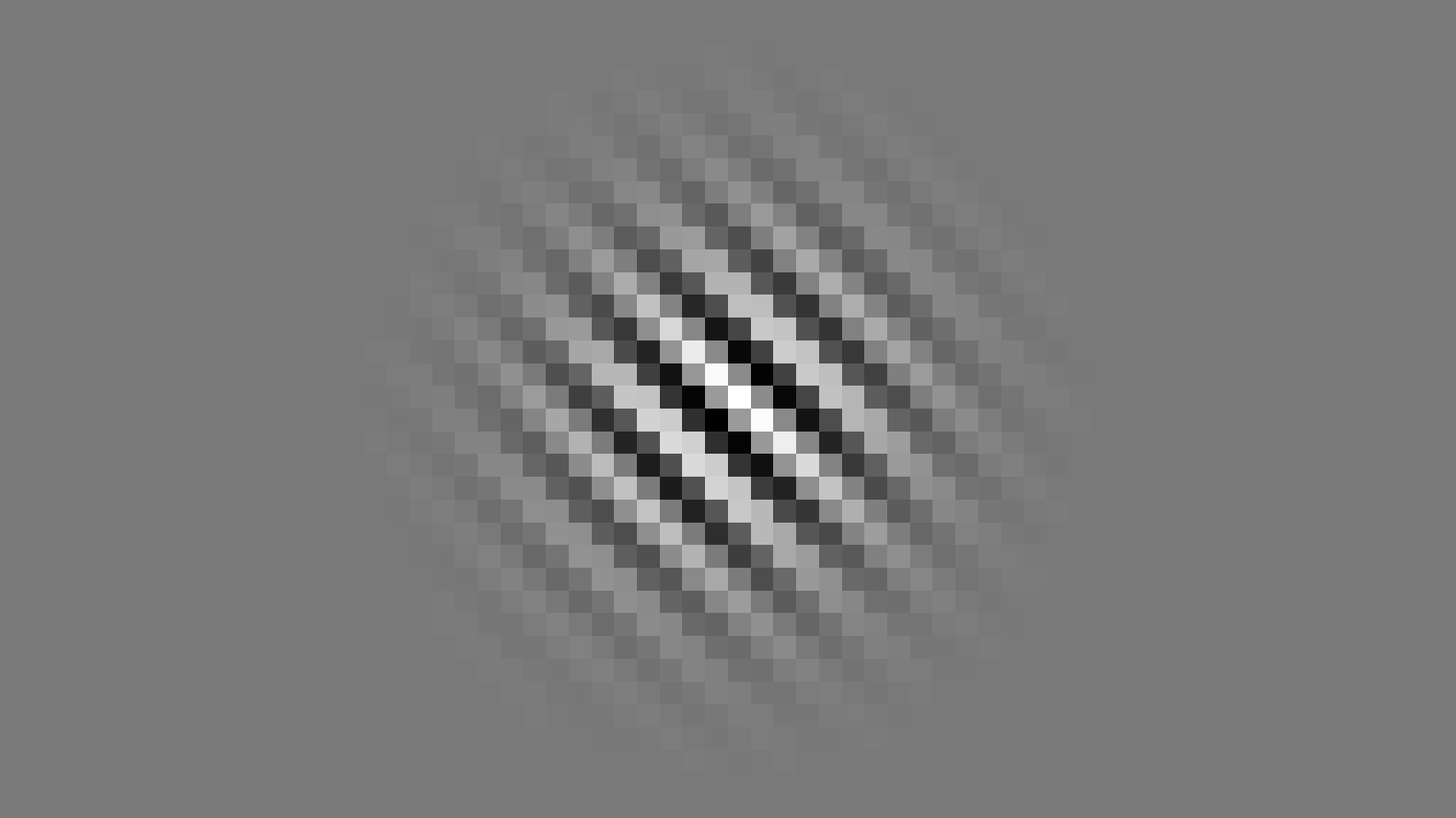}
    \includegraphics[width=0.45\textwidth]{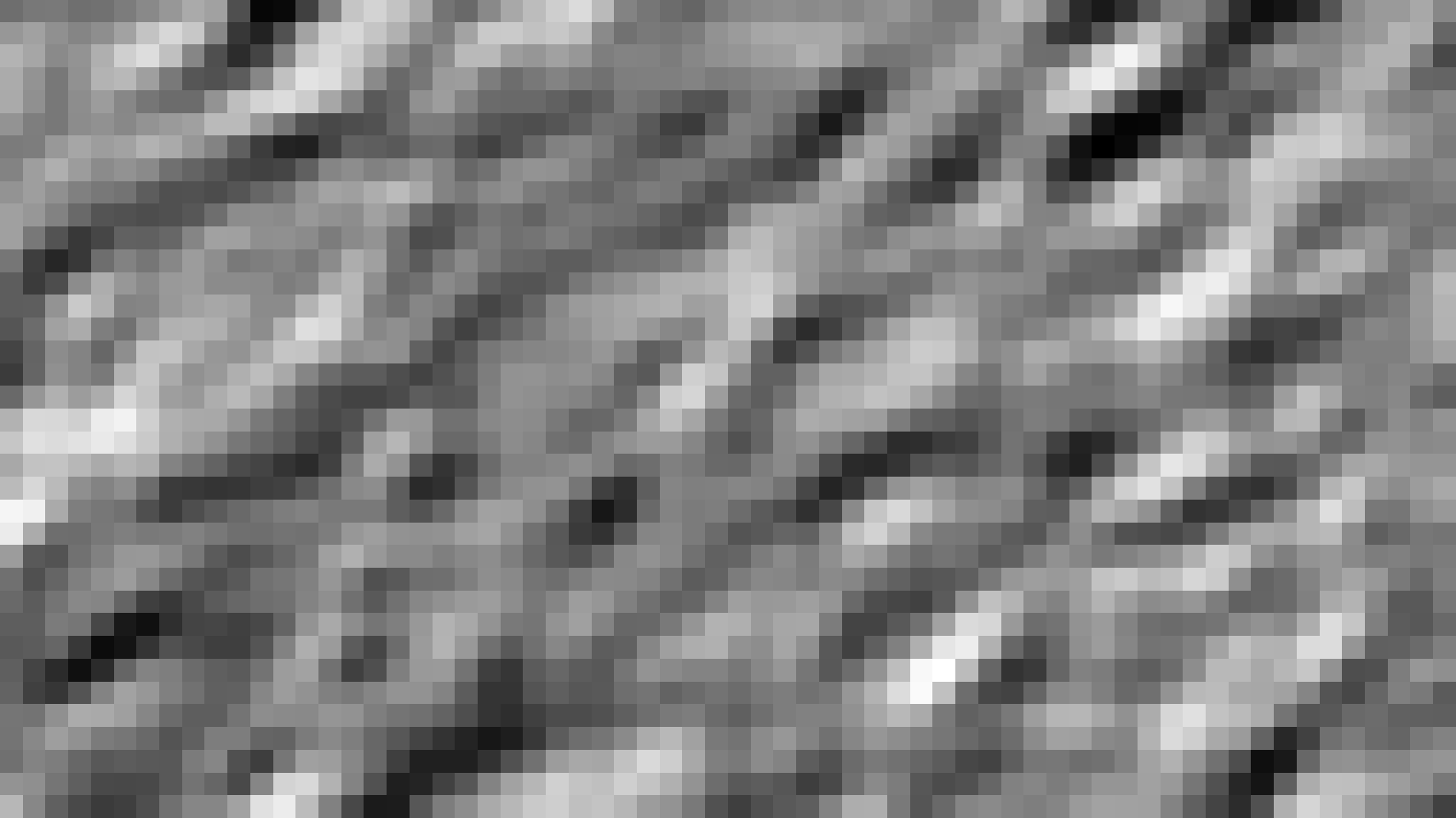}
    \includegraphics[width=0.45\textwidth]{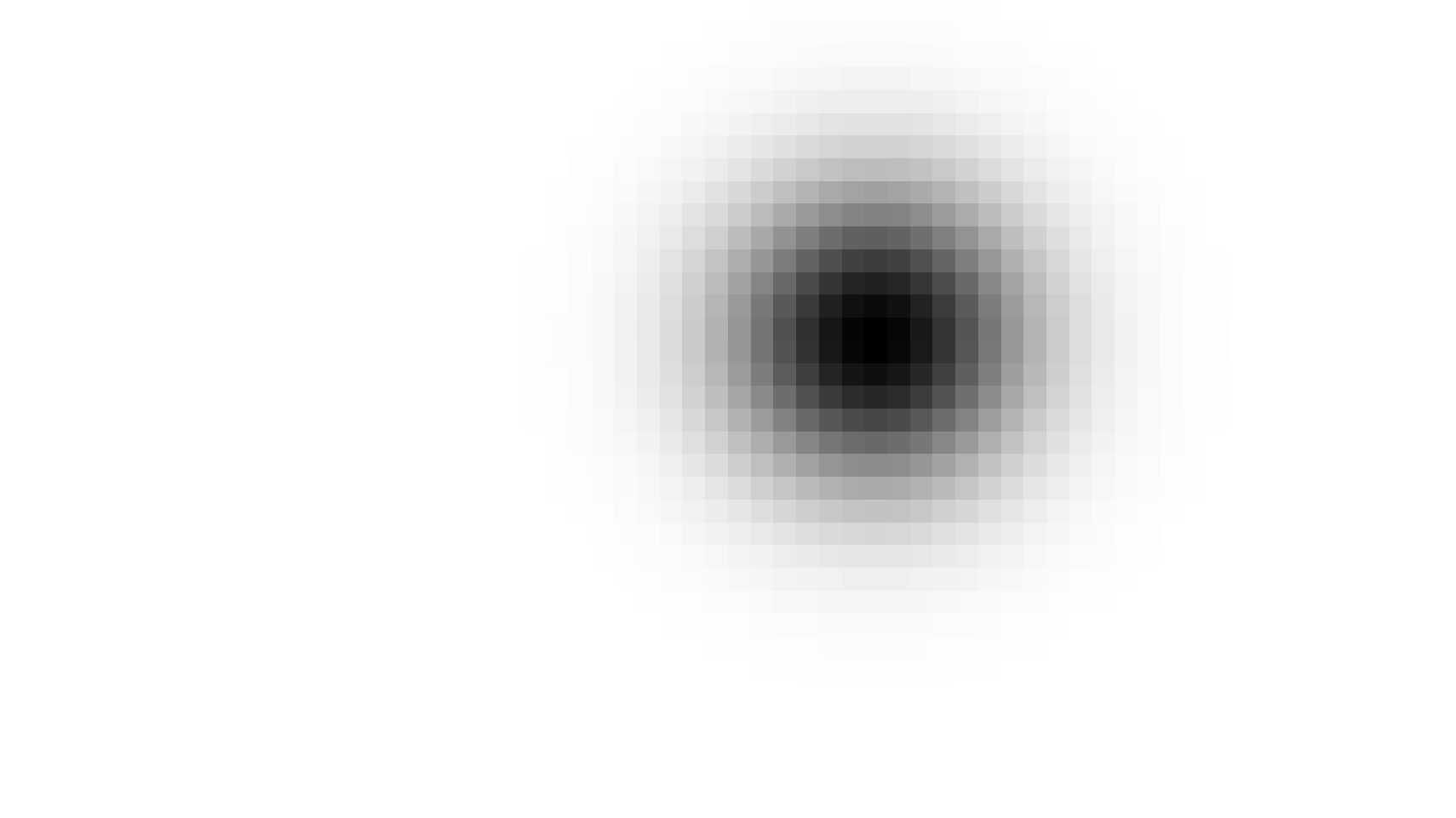}
    \caption{Middle frame for each video type considered in this work: \textsc{Naturalistic} (upper left), \textsc{Gaussian} (upper right), \textsc{Waves} (bottom left), \textsc{Moving Dot} (bottom right).}
    \label{fig:videotypes}
\end{figure}

In this section, we provide further details on the Sensorium dataset \cite{sensorium2023} used in this work.
\paragraph{Video classes.} For each mouse, $\approx 700$ videos are provided with corresponding neural responses and behavioural metadata. These videos are characterized by subgroups of repeated entries. Video types change across mice, we select $4$ video types and label them as \textsc{Naturalistic}, \textsc{Gaussian}, \textsc{Waves} and \textsc{Moving Dot}. We show a sample frame for each type in Figure \ref{fig:videotypes}. Groups of repeated videos have typically $10$ samples each. Each video contains 324 frames, with the last 24 frames being blank. The format of each video is $(36,64)$ in grayscale.

\begin{figure}[h!]
    \centering
    \includegraphics[width=0.95\textwidth]{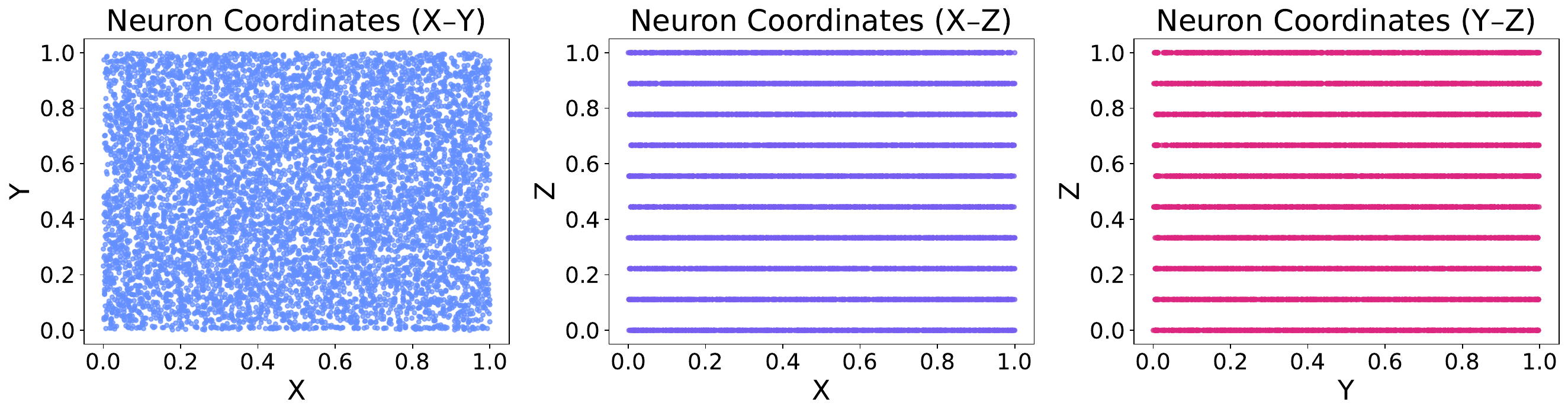}
    \caption{
\textbf{Spatial sampling geometry.} 
Neuron positions for a representative mouse–video pair, shown in three projections: X–Y (left), X–Z (center), and Y–Z (right). The laminar stack of regularly spaced z-planes is visible, and within-plane coverage appears uniform. Coordinates are normalized to the imaging volume.
}
\label{fig:neuron_coord}
\end{figure}

\begin{figure}
    \centering
    \includegraphics[width=0.45\textwidth]{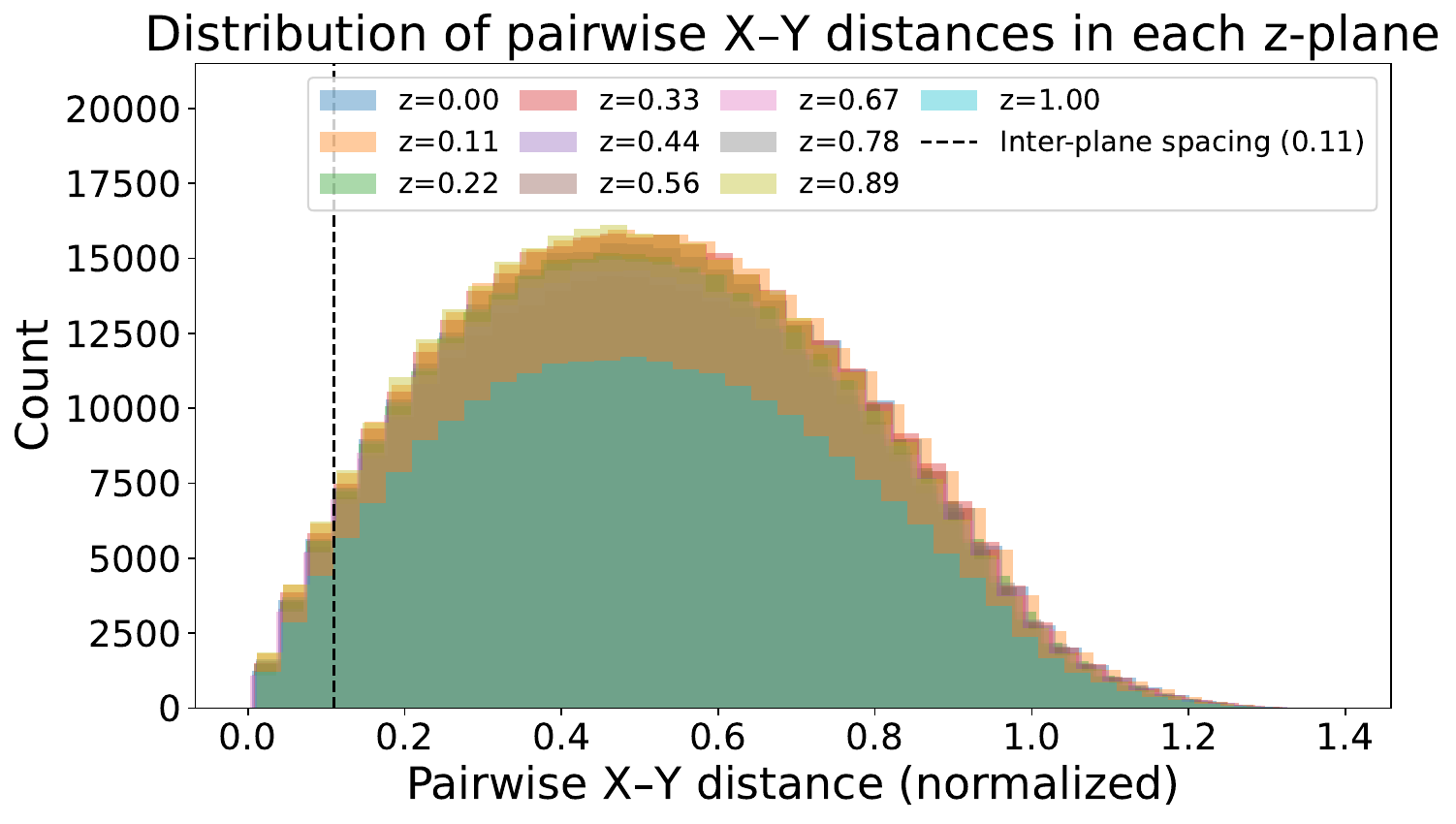}
    \includegraphics[width=0.45\textwidth]{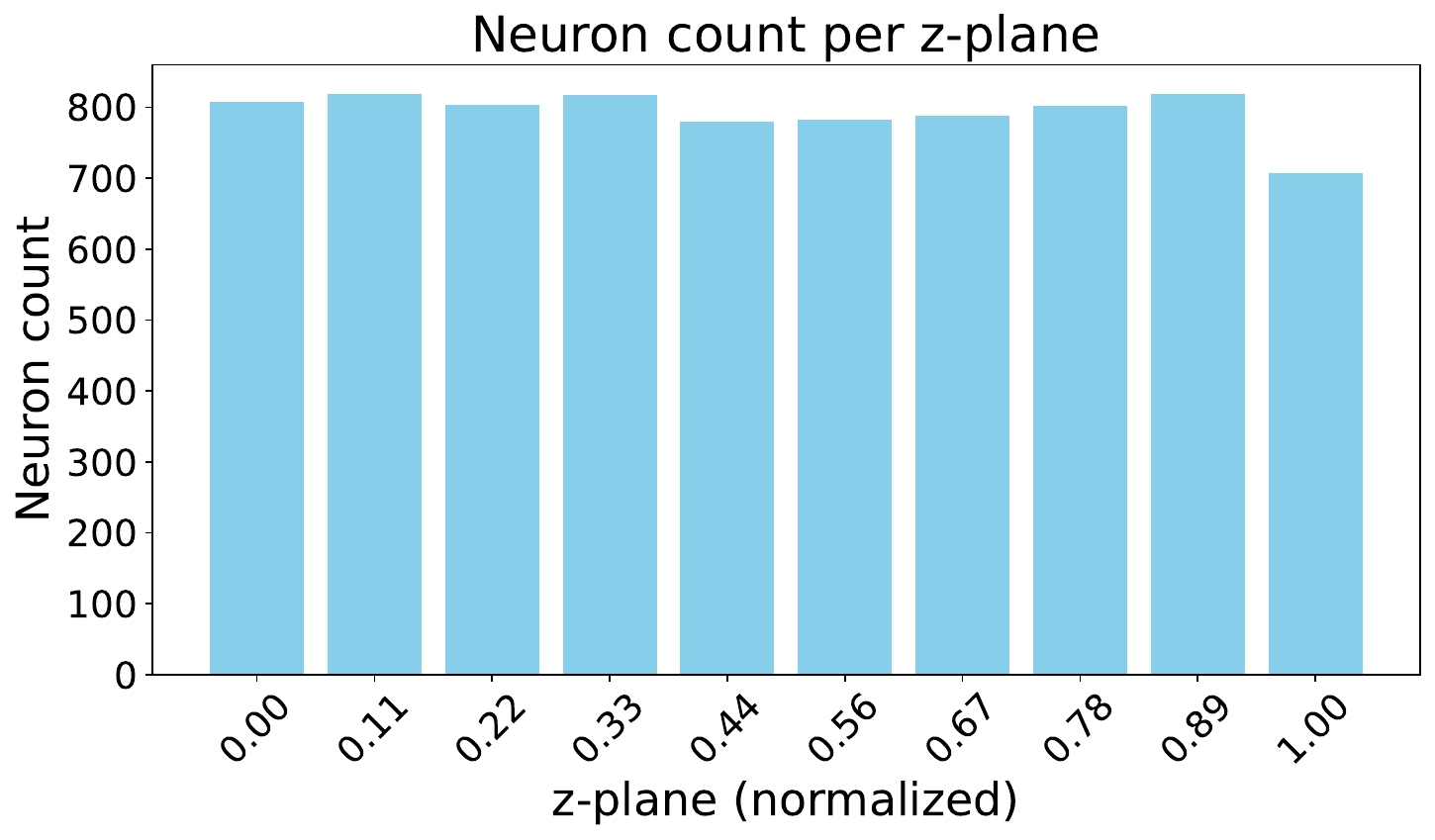}
    \caption{Left:
    \textbf{Within-plane sampling density.}
    Distribution of pairwise X–Y distances between neurons for each z-plane (overlaid histograms). The vertical dashed line marks the normalized inter-plane spacing ($\approx 0.111$), serving as a reference scale. The similar distributions across planes indicate uniform sampling density within each image plane. Right:
    \textbf{Depth sampling balance.}
Bar plot of neuron counts per z-plane for the same mouse–video pair. Cell counts are broadly balanced across planes, with only minor variability, confirming uniform sampling in depth.
}
\label{fig:counts}
\end{figure}

\paragraph{Spatial sampling and geometry of neuron positions.} For a representative mouse–video pair (mouse 2-10, video 34), we show cell positions in the three projections (X–Y, X–Z, Y–Z) in Figure \ref{fig:neuron_coord}.
Neurons populate a stack of regularly spaced z-planes; after coordinate normalization, the spacing between adjacent planes is $\approx 0.11$ (marked in the distance plots as a reference scale). Within-plane sampling density appears comparable across slices: pairwise X–Y distances have similar distributions for each z-plane, indicating a fairly uniform coverage in each image plane. Cell counts per plane are also broadly balanced, with only minor variability across z. We show these characteristics for the same sample case in Figure \ref{fig:counts}.
We also checked for potential duplicates across adjacent planes: only two (x, y) locations recur at different z, and their fluorescence time series are numerically indistinguishable, indicating rare double detections that are negligible at the dataset scale.

\begin{figure}[h!]
    \centering
    \includegraphics[width=0.45\textwidth]{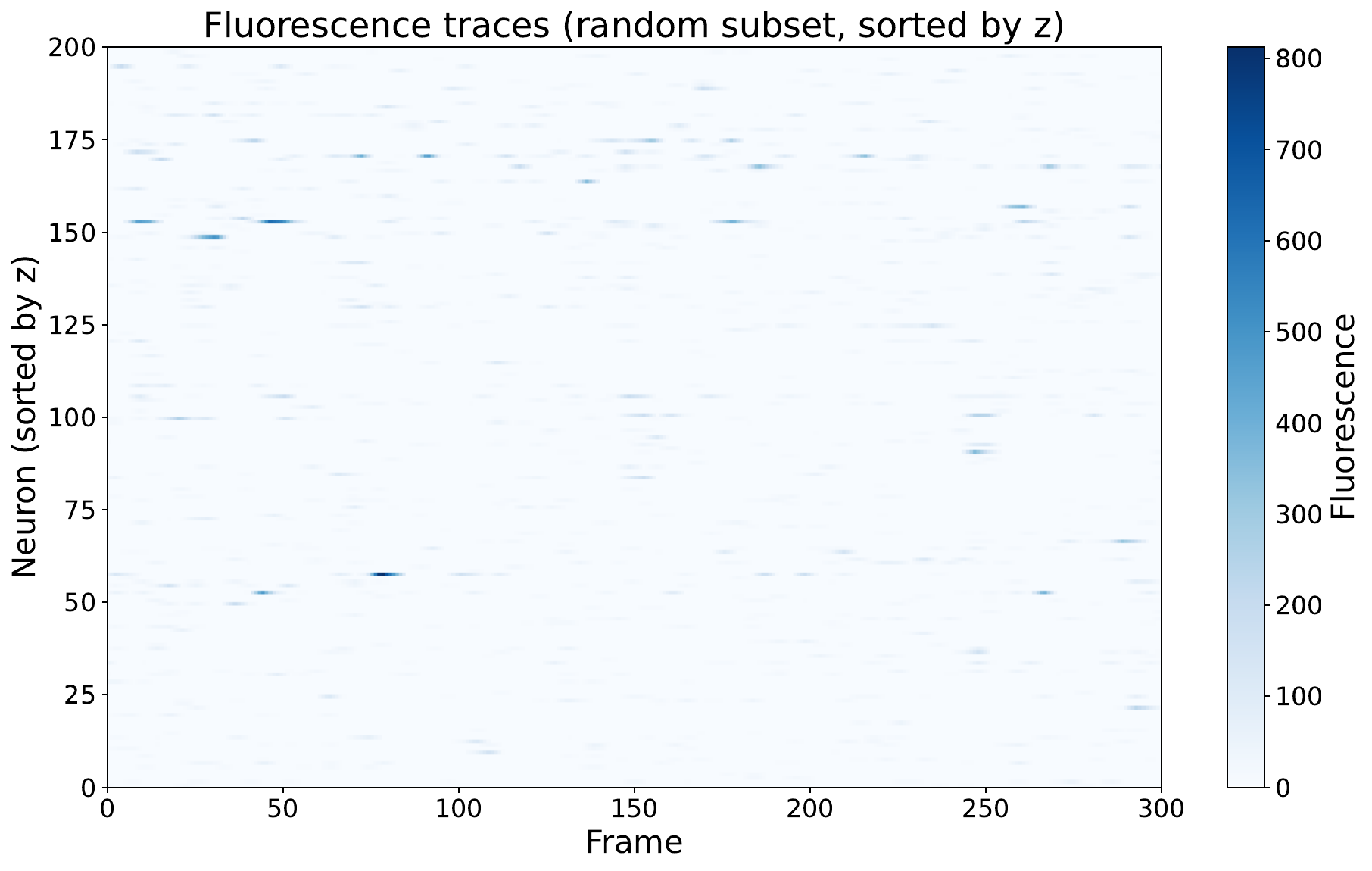}
    \includegraphics[width=0.45\textwidth]{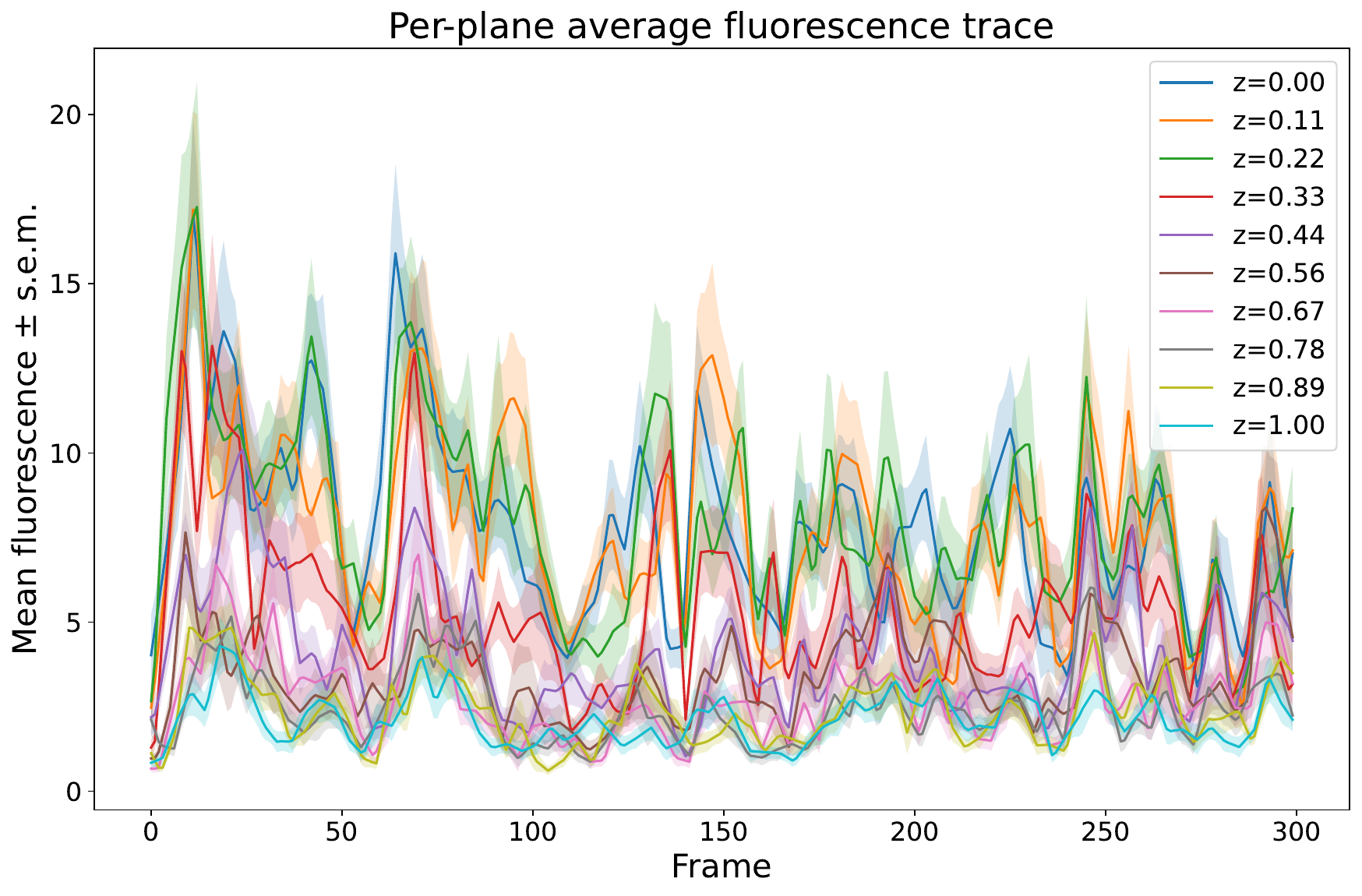}
    \caption{Left:
    \textbf{Fluorescence response overview.}
Heatmap of single-trial fluorescence traces for a random subset of 200 neurons, sorted by z-plane. Each row is a neuron, each column a video frame. This illustrates the diversity and temporal structure of neural responses across the imaging volume. Right:
\textbf{Per-plane average fluorescence dynamics.}
Mean fluorescence trace (solid line) and standard error (shaded area) for each z-plane, showing depth-dependent response dynamics and variability across neurons and planes.
}
\label{fig:traces}
\end{figure}

\paragraph{Neural responses.} For each neuron we use single-trial fluorescence traces time-locked to the video frames. The duplicate check above confirms that when a cell is inadvertently captured in neighboring planes, its activity traces are identical, as expected. In the main analyses we work with per-frame response interpolated on 2D grids; for descriptive purposes, it is useful to visualize raw traces and simple aggregates per plane to illustrate trial structure and variability across neurons and depths, see Figure \ref{fig:traces}.

\section{Zigzag algorithm}\label{app:algo}
\begin{figure}[h!]
\includegraphics[width=0.99\textwidth]{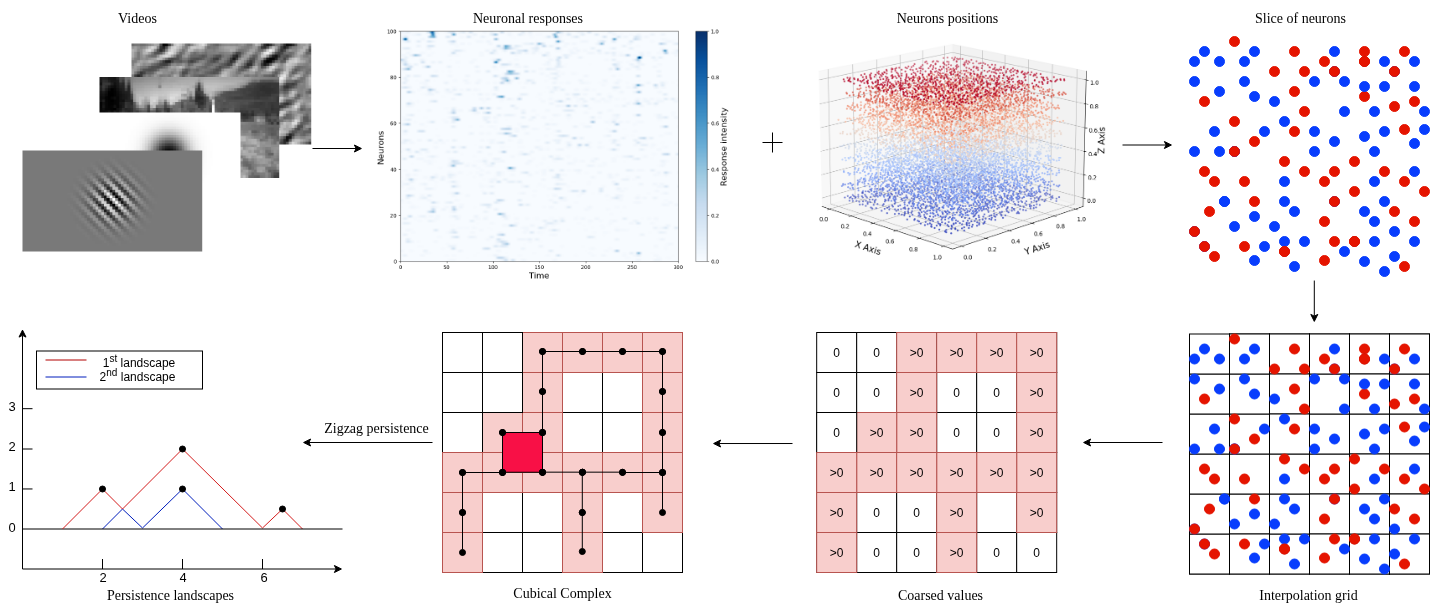}
\caption{Schematic representation of the full zigzag pipeline.}
\label{fig:pipeline}
\end{figure}

Here we describe the full analysis pipeline. A schematic representation of each step is shown in Figure \ref{fig:pipeline}.
The algorithm is presented below:
\begin{algorithm}[h!]
\caption{Zigzag persistence pipeline for frame-wise neuronal activity grids}
\begin{algorithmic}[1]
\Require List of videos \texttt{all\_indices}; number of z-planes $Z$; grid size $n_{\rm grid}$; per-video, per-plane interpolated activity grids of shape $(n_{\rm grid}, n_{\rm grid}, T)$; landscape params $(R,k)$.
\For{each video index $i$ in \texttt{all\_indices}}
\For{each z-plane $z = 0, \dots, Z-1$}
\State \texttt{path} $\gets$ path to interpolated activity grid for $(i,z)$
\If{\texttt{path} does not exist} \State \textbf{continue} \EndIf
\State $r \gets$ load activity grid $(n_{\rm grid}, n_{\rm grid}, T)$
\State $\delta \gets \mathrm{normalization}(r)$ \Comment{$\delta_{a,z,i}(T) = (r_{a,z,i}(T) - \bar{r}_{a,z,i})/\bar{r}_{a,z,i}$}
\State Optional control: frame shuffle $t \mapsto \pi(t)$; grid scramble via fixed spatial permutation $p$ per $(i,z)$
\State \texttt{complexes} $\gets$ [ ]
\For{$t = 1, \dots, T$}
\State \texttt{mask} $\gets [\,\delta(T) > 0\,]$
\State Build cubical cells active at $T$: vertices, grid-edges, unit-squares from \texttt{mask}
\State Build simplicial complex $\mathcal S_t$ by an abstract closure adapter:
\State \quad insert all active vertices as 0-simplices
\State \quad insert all active grid-edges as 1-simplices
\State \quad for each active square with corners ${v_{00},v_{10},v_{01},v_{11}}$:
\State \qquad insert the abstract 3-simplex on these four vertices and retain only its 2-skeleton
\State append $\mathcal S_t$ to \texttt{complexes}
\EndFor
\State Build interleaved zigzag sequence:
\State \quad $X \gets [\,\mathcal{S}_1,\ \mathcal{S}_1 \cap \mathcal{S}_2,\ \mathcal{S}_2,\ \dots,\ \mathcal{S}_{T-1}\cap \mathcal{S}_T,\ \mathcal{S}_T\,]$
\State $f,\ \texttt{times} \gets \mathrm{encodeZigzagFiltration}(X)$
\State \texttt{bars} $\gets \mathrm{FastZigZag}(f,\ \texttt{times})$
\State \texttt{bars\_$H_1$} $\gets \mathrm{selectDimension}(\texttt{bars},,1)$
\State Vectorize \texttt{bars\_$H_1$} as persistence landscapes: sample $R$ points, keep $k$ layers
\State Save per-plane descriptor for $(i,z)$
\EndFor
\EndFor
\end{algorithmic}
\end{algorithm}

For each mouse $a$, video $v$, plane $p$ and frame $t$, we load the interpolated activity grid $r_{a,p,v}(t)$ and convert it to a per-pixel activity
$\delta_{a,p,v}(t)$ which centers each pixel at zero and makes the threshold $\delta > 0$ correspond to above-average activity within that plane and video.

To probe what the descriptor uses, we optionally apply:
\begin{itemize}
\item a frame-order shuffle $t \mapsto \pi(t)$ (identical across planes) to disrupt temporal order while preserving per-frame spatial structure;
\item a fixed spatial permutation $\pi$ per $(i,z)$ to scramble the grid indices across all frames, disrupting spatial contiguity while largely preserving marginal activity statistics.
\end{itemize}

\noindent\textbf{Per-frame complexes and simplicial adapter.}
For each frame $t$, we form a binary mask of active pixels ($\delta > 0$) and construct the corresponding 2D cubical cells (active vertices, grid-edges, and unit-squares). To interface with zigzag software operating on simplicial complexes, we map this cubical complex to a simplicial complex $S_t$ using an abstract closure:
(i) insert all active vertices as 0-simplices and all active grid-edges as 1-simplices; (ii) for each active square with four grid-corner vertices, insert the abstract 3-simplex on those four vertices and retain only its 2-skeleton (thus adding both diagonals and all triangular faces). This is not a triangulation; it over-connects to ensure consistent 2D faces across frames and stable computation. Vertex identifiers are the grid coordinates, so identical simplices are recognized across time.

\noindent\textbf{Interleaved intersection zigzag.}
We assemble the time sequence with intersections at intermediate steps:
\begin{equation}
\mathcal{S}_1 \ \hookleftarrow\ \mathcal{S}_1 \cap \mathcal{S}_2\ \hookrightarrow\ \mathcal{S}_2 \ \hookleftarrow\ \mathcal{S}_2 \cap \mathcal{S}_3\ \hookrightarrow\ \cdots \ \hookleftarrow\ \mathcal{S}_{T-1}\cap \mathcal{S}_T\ \hookrightarrow\ \mathcal{S}_T.
\end{equation}
Intersections are set-theoretic intersections of simplices using global vertex labels. This produces a valid zigzag of inclusions that captures gains and losses of simplices between consecutive frames while maintaining well-defined maps.

\noindent\textbf{Zigzag persistence and vectorization.}
We compute $H_1$ zigzag persistent homology on the sequence $X$ using two public libraries:
\begin{itemize}
\item \textsc{Dionysus2} \cite{dionysus}, which provides data structures and routines to encode the zigzag filtration and its time indices (birth/death layers);
\item \textsc{fastzigzag} \cite{dey2022fastcomputationzigzagpersistence}, which converts the zigzag filtration to an equivalent non-zigzag filtration of an auxiliary complex, computes standard persistence efficiently, and maps intervals back to the zigzag setting to obtain the $H_1$ barcode.
\end{itemize}
This combination of codes and algorithm allows for a computational time of about $150$ ms on a single CPU per single data vector, i.e. roughly $100$ minutes to run the full pipeline on a single mouse ($\approx 700$ video stimuli).

\begin{figure}
\centering
    \includegraphics[width=0.45\textwidth]{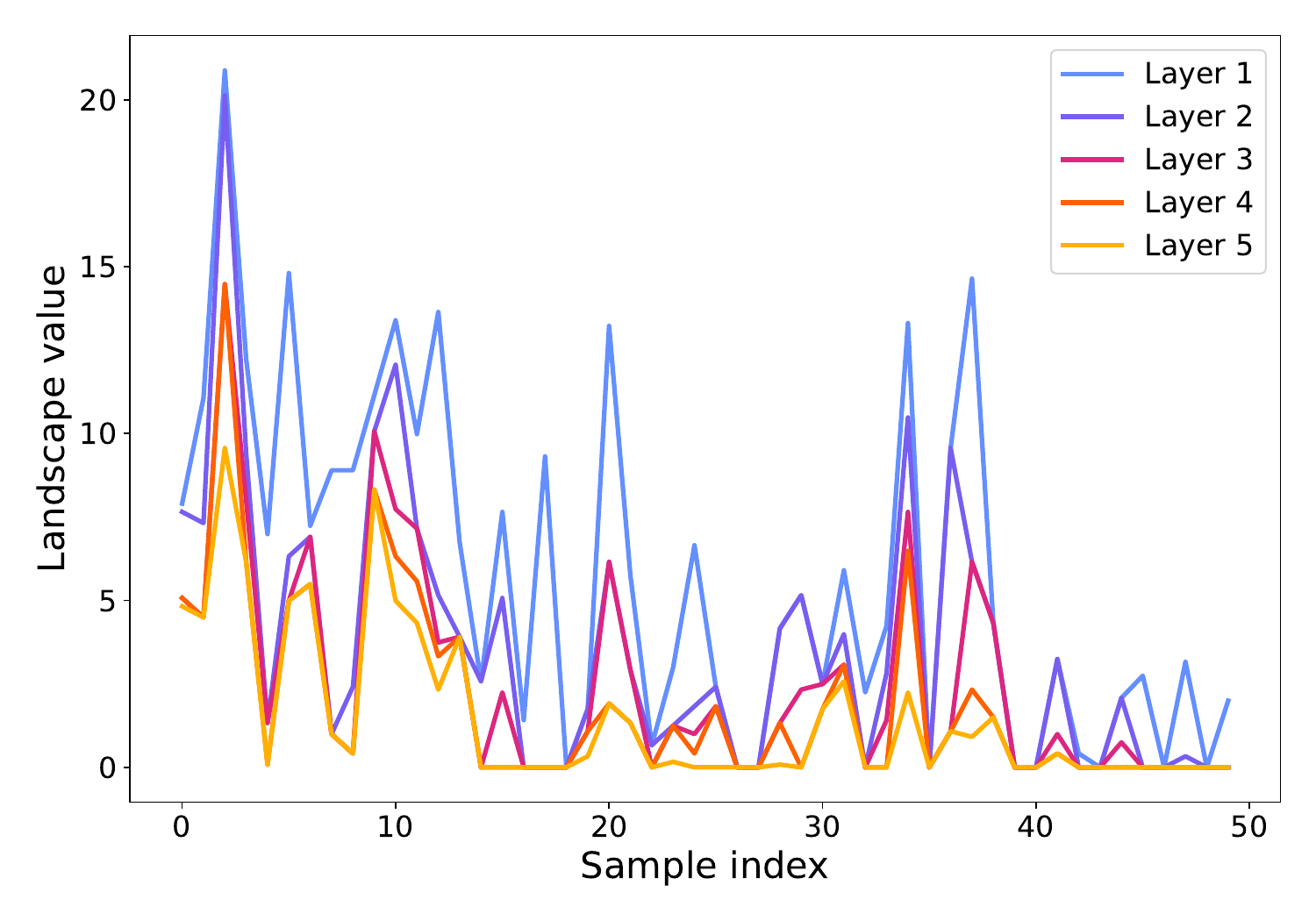}
    \includegraphics[width=0.45\textwidth]{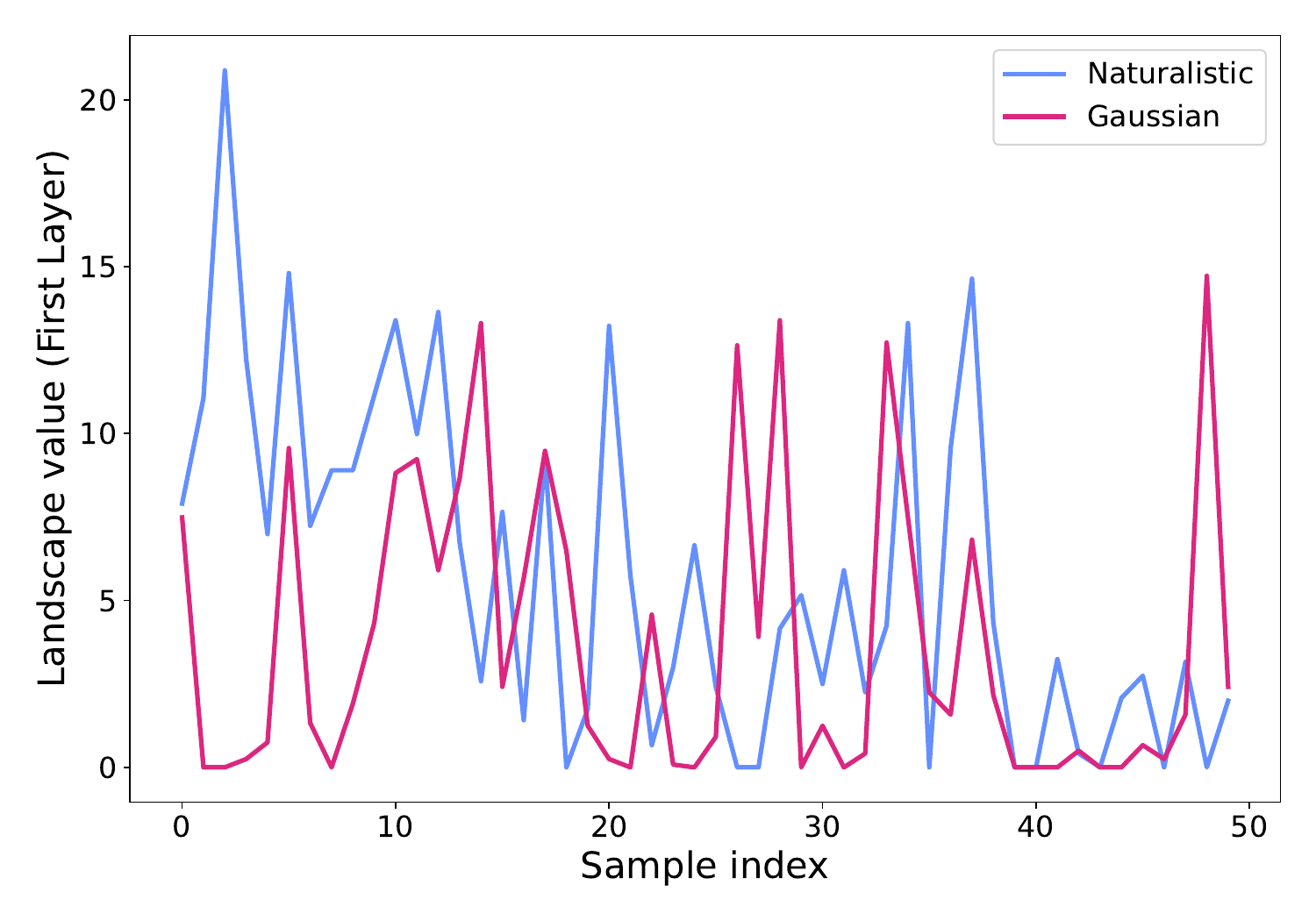}
    \caption{\textbf{Representative persistence landscape vectors.}\\
Left: Each curve shows one of the first five landscape layers ($\Lambda_k$) for a single trial and a single $z$-plane, summarizing the time-varying prominence of topological cycles ($H_1$) in the neuronal activity. Right: Each curve shows the first landscape layer ($\Lambda_1$) of a single $z$-plane of trials corresponding to a naturalistic video and a Gaussian video. All trials refer to mouse \textsc{2-10}. }
\label{fig:landscape}
\end{figure}

\paragraph{Topological descriptors.} Each $H_1$ barcode is mapped to persistence landscapes $\Lambda_k$ sampled on a uniform grid over the time index; in our experiments we use $R=50$ sample points and $k=5$ landscape layers, yielding a 250-dimensional vector per plane. Trial-level descriptors concatenate all planes. In Figure \ref{fig:landscape} we show a representative data vector for mouse $\textsc{2-10}$ on a single $z$-plane both for the full 250-dimensional vector for a naturalistic video and a comparison of the first landscapes between a naturalistic video and a Gaussian video.

\section{Additional results}\label{app:results}

In this section, we quote additional results to complement the results quoted in the main text. 

\begin{figure}
\centering
\includegraphics[width=0.4\textwidth]{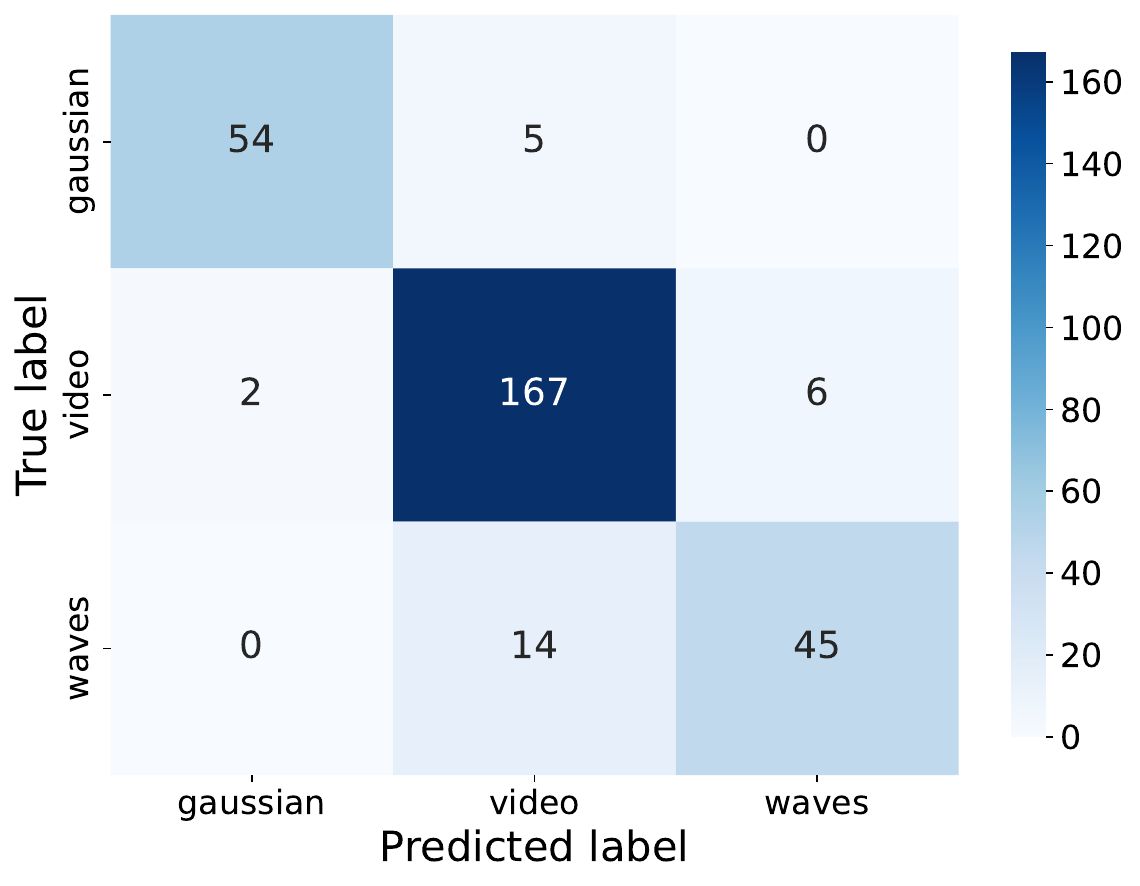}
\includegraphics[width=0.4\textwidth]{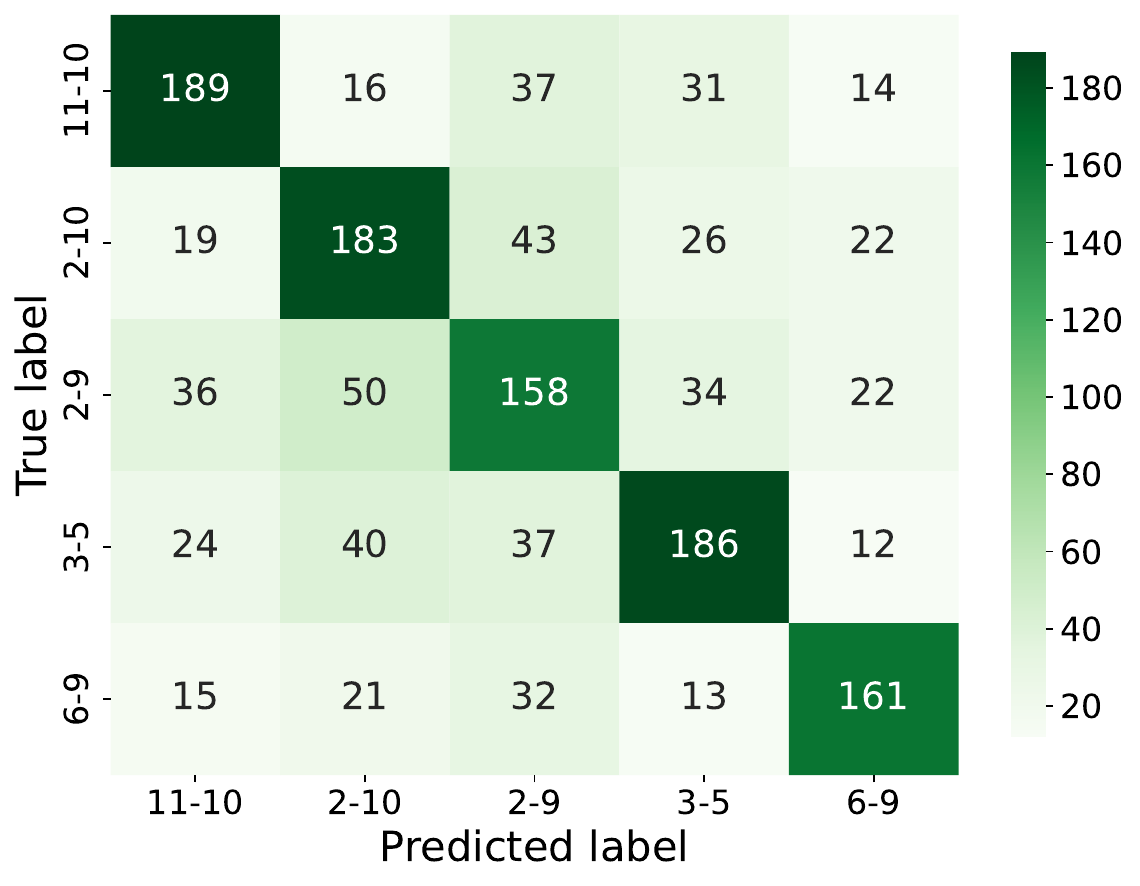}
\caption{Aggregated confusion matrix for 5 train/test splits. \emph{Left:} video type. \emph{Right:} Mouse identity.}
\label{fig:classifier}
\end{figure}

\paragraph{Clustering protocol (A).} For the clustering experiments, we run the pipeline on all $5$ mice, shown in Table \ref{tab:allmice}, where we also include the accuracy of classification by clustering as a metric. For these, we did not run control-type checks. Results are broadly coherent with what shown in Table \ref{tab:exp_ari_main} for mouse \textsc{2-10}.
\begin{table}[h]
\centering
\small
\begin{tabular}{lcc}
\toprule
Mouse \textsf{2-10} & & \\
\toprule
Condition & Baseline ARI & Baseline Accuracy \\
\midrule
  Naturalistic & 0.945 $\pm$ 0.118 & 0.955 $\pm$ 0.076  \\
  Gaussian & 0.702 $\pm$ 0.156 & 0.837 $\pm$ 0.104\\
  Waves & 0.603 $\pm$ 0.218 & 0.799 $\pm$ 0.126  \\
\bottomrule
\toprule
Mouse \textsf{2-9} & & \\
\toprule
Condition & Baseline ARI & Baseline Accuracy \\
\midrule
  Naturalistic & 0.734 $\pm$ 0.161 & 0.863 $\pm$ 0.103  \\
  Gaussian & 0.820 $\pm$ 0.141 & 0.921 $\pm$ 0.071\\
  Moving Dot & 0.463 $\pm$ 0.239 & 0.720 $\pm$ 0.148  \\
\bottomrule
\toprule
Mouse \textsf{11-10} & & \\
\toprule
Condition & Baseline ARI & Baseline Accuracy \\
\midrule
  Naturalistic & 0.791 $\pm$ 0.185 & 0.884 $\pm$ 0.117  \\
  Waves & 0.423 $\pm$ 0.225 & 0.683 $\pm$ 0.128  \\
  Moving Dot & 0.982 $\pm$ 0.069 & 0.993 $\pm$ 0.026  \\
\bottomrule
\toprule
Mouse \textsf{3-5} & &\\
\bottomrule
Condition & Baseline ARI & Baseline Accuracy \\
\midrule
  Naturalistic & 0.911 $\pm$ 0.098 & 0.963 $\pm$ 0.044  \\
  Waves & 0.446 $\pm$ 0.196 & 0.706 $\pm$ 0.120  \\
  Moving Dot & 0.791 $\pm$ 0.078 & 0.921 $\pm$ 0.031  \\
\bottomrule
\toprule
Mouse \textsf{6-9} & & \\
\bottomrule
Condition & Baseline ARI & Baseline Accuracy \\
\midrule
  Naturalistic & 0.912 $\pm$ 0.138 & 0.955 $\pm$ 0.081  \\
  Gaussian & 0.953 $\pm$ 0.067 & 0.983 $\pm$ 0.025\\
\bottomrule
\end{tabular}
\caption{Clustering ARI and accuracy (mean $\pm$ std over 20 runs) for the all mice.}
\label{tab:allmice}
\end{table}
\paragraph{Classification protocol (B) video family.} Here we provide for all mice the F1-score for each class and the cross-validation accuracies on the task of predicting the right video family, Table \ref{tab:per_mouse_f1}. Additionally, we show the aggregated confusion matrix for 5 train/test splits for the mouse \textsc{2-10}, see left panel of Figure \ref{fig:classifier}.\footnote{Note: discrepancy of accuracy rates between the aggregated confusion matrix and the cross-validation score is because the latter uses mean per-fold accuracy, which is more reliable for unbalanced classes and small datasets, as in this case. The aggregated confusion matrix is quoted for interpretability purposes.}
\begin{table}[h]
\centering
\small
\begin{tabular}{l c c}
\toprule
Mouse \textsc{2-10} & & \\
\midrule
 \textbf{Class} & \textbf{F1} & \textbf{Support} \\
\midrule
Naturalistic & 0.88 & 35 \\
Gaussian    & 0.82 & 12 \\
Waves    & 0.78 & 12 \\
\cmidrule(lr){1-3}
\multicolumn{3}{r}{\textbf{CV accuracy:} 0.686 $\pm$ 0.055} \\
\bottomrule
\toprule
Mouse \textsc{2-9} & & \\
\midrule
 \textbf{Class} & \textbf{F1} & \textbf{Support} \\
\midrule
Naturalistic & 0.95 & 37 \\
Gaussian    & 0.96 & 12 \\
Moving Dot    & 0.83 & 11 \\
\cmidrule(lr){1-3}
\multicolumn{3}{r}{\textbf{CV accuracy:} 0.743 $\pm$ 0.064} \\
\bottomrule
\toprule
Mouse \textsc{11-10} & & \\
\midrule
 \textbf{Class} & \textbf{F1} & \textbf{Support} \\
\midrule
Naturalistic & 0.93 & 35 \\
Gaussian    & 0.86 & 11 \\
Moving Dot    & 0.91 & 12 \\
\cmidrule(lr){1-3}
\multicolumn{3}{r}{\textbf{CV accuracy:} 0.721 $\pm$ 0.115} \\
\bottomrule
\toprule
Mouse \textsc{3-5} & & \\
\midrule
 \textbf{Class} & \textbf{F1} & \textbf{Support} \\
\midrule
Naturalistic & 0.89 & 36 \\
Waves    & 0.73 & 12 \\
Moving Dot    & 0.82 & 12 \\
\cmidrule(lr){1-3}
\multicolumn{3}{r}{\textbf{CV accuracy:} 0.672 $\pm$ 0.056} \\
\bottomrule
\toprule
Mouse \textsc{6-9} & & \\
\midrule
 \textbf{Class} & \textbf{F1} & \textbf{Support} \\
\midrule
Naturalistic & 1.0 & 37 \\
Gaussian    & 1.0 & 12 \\
\cmidrule(lr){1-3}
\multicolumn{3}{r}{\textbf{CV accuracy:} 0.888 $\pm$ 0.064} \\
\bottomrule
\end{tabular}
\caption{Per-class F1 scores and supports per mouse for a single train/test split. Cross-validated accuracy is run on 5 train/test splits.}
\label{tab:per_mouse_f1}
\end{table}

\paragraph{Classification protocol (C) mouse identity.} Here we provide for the aggregated confusion matrix for 5 train/test splits on the mouse identity classification task, see Right panel of Figure \ref{fig:classifier}.

\end{document}